\documentclass[12pt,reqno]{amsart}
\usepackage{latexsym,amsmath,amsfonts,amssymb,amsthm}
\theoremstyle{plain}
\newtheorem{Thm}{Theorem}
\newtheorem{Lem}{Lemma}
\newtheorem{Cor}{Corollary}

\theoremstyle{definition}

\theoremstyle{remark}

\def\Z{\mathbb Z}
\def\N{\mathbb N}

\def\P{\mathcal P}
\def\pmod #1{\ ({\rm mod}\ #1)}

\def\floor #1{\lfloor{#1}\rfloor}
\def\ceil #1{\left\lceil{#1}\right\rceil}
\def\1{\mathbf 1}

\begin{document}
\title{A density version of Vinogradov's three primes theorem}
\author{Hongze Li}
\email{lihz@sjtu.edu.cn}
\address{
Department of Mathematics, Shanghai Jiaotong University, Shanghai
200240, People's Republic of China}
\author{Hao Pan}
\email{haopan79@yahoo.com.cn}
\address{
Department of Mathematics, Nanjing University, Nanjing 210093,
People's Republic of China}
 \thanks{This work was supported by the
National Natural Science Foundation of China (Grant No. 10771135).
The second author is the corresponding author.}\maketitle
\begin{abstract}
Let $\P$ denote the set of all primes. Suppose that $P_1$, $P_2$,
$P_3$ are three subsets of $\P$ with $\underline{d}_{\P}(P_1)+
\underline{d}_{\P}(P_2)+\underline{d}_{\P}(P_3)>2$, where
$\underline{d}_{\P}(P_i)$ is the lower density of $P_i$ relative to
$\P$. We prove that for every sufficiently large odd integer $n$,
there exist $p_i\in P_i$ such that $n=p_1+p_2+p_3$.
\end{abstract}
\section{Introduction}
\setcounter{equation}{0} \setcounter{Thm}{0} \setcounter{Lem}{0}
\setcounter{Cor}{0}

The ternary Goldbach conjecture says that every odd integer
greater than 7 is the sum of three primes. This problem was
basically solved by Vinogradov \cite{Vinogradov37} in 1937, and in
fact he showed that for every sufficiently large odd integer $n$,
$$
\sum_{\substack{p_1+p_2+p_3=n\\ p_1, p_2, p_3 \text{ prime}}}\log
p_1\log p_2\log p_3=\frac{1}{2}\mathfrak{S}(n)n^2+O(n^2(\log
n)^{-A}),
$$
where
$$
\mathfrak{S}(n)=\prod_{p\nmid n}(1+(p-1)^{-3}) \prod_{p\mid
n}(1-(p-1)^{-2})
$$
and $A$ is a positive constant. Nowadays Vinogradov's theorem has
become a classical result in additive number theory. Later, using
a similar method, van der Corput \cite{Corput39} proved that the
primes contain infinitely many non-trivial 3-term arithmetic
progressions (3AP).

On the other hand, another classical result due to Roth
\cite{Roth53} asserts that any subset $A$ of the integers with
$\overline{d}(A)>0$ contains infinitely many non-trivial 3APs, where
$$
\overline{d}(A)=\limsup_{x\to\infty}\frac{|A\cap[1,x]|}{x}.
$$
Roth's theorem is a special case of the well-known Szemer\'edi
theorem \cite{Szemeredi75}, which states that any integers set $A$
with $\overline{d}(A)>0$ contains arbitrarily long arithmetic
progressions.

For two non-empty sets $X$ and $A$ of positive integers, define
the upper density and lower density of $A$ relative to $X$ by
$$
\overline{d}_X(A)=\limsup_{x\to\infty}\frac{|A\cap
X\cap[1,x]|}{|X\cap[1,x]|}
$$
and
$$
\underline{d}_X(A)=\liminf_{x\to\infty}\frac{|A\cap
X\cap[1,x]|}{|X\cap[1,x]|}.
$$
Let $\P$ denote the set of all primes. In \cite{Green05}, Green
obtained a Roth-type generalization of van der Corput's result.
Green showed that if $\P_0$ is a subset of $\P$ with
$\overline{d}_{\P}(\P_0)>0$ then $\P_0$ contains infinitely many
3APs. One major ingredient in Green's proof is a transference
principle, which transfers a subset of primes with relative positive
density to a subset of $\Z_N=\Z/N\Z$ (where $N$ is a large prime)
with positive density. Subsequently, this principle was greatly
improved (in a different way) in the proof of Green and Tao's
celebrated theorem \cite{GreenTao} that the primes contain
arbitrarily long arithmetic progressions.

The Hardy-Littlewood circle method \cite{Vaughan97} is commonly
applied in Vinogradov's, van der Corput's, Roth's and Green's
proofs. In this paper, we shall use Green's idea to extend the
Vinogradov theorem as follows.
\begin{Thm}
\label{t1} Suppose that $P_1, P_2, P_3$ are three subsets of $\P$
with
$$
\underline{d}_{\P}(P_1)+
\underline{d}_{\P}(P_2)+\underline{d}_{\P}(P_3)>2.
$$
Then for every sufficiently large odd integer $n$, there exist
$p_1\in P_1$, $p_2\in P_2$ and $p_3\in P_3$ such that
$n=p_1+p_2+p_3$.
\end{Thm}

Notice that the result of Theorem \ref{t1} is the best possible in
the following sense:

Letting $P_1=P_2=\{p\in\P:\,p\equiv1\pmod{3}\}$ and
$P_3=\P\setminus\{3\}$, then
$\underline{d}_{\P}(P_1)=\underline{d}_{\P}(P_2)=1/2$ and
$\underline{d}_{\P}(P_3)=1$, but $6k+5\not\in P_1+P_2+P_3$ for any
integer $k$.

For a positive integer $q$, let $\Z_q=\Z/q\Z$ and
$\Z_q^*=\{b\in\Z_q:\, (b,q)=1\}$. The key of our proof is an
addition theorem:
\begin{Thm}
\label{t2} Let $q$ be a positive integer with $(q,6)=1$. Let $f_1,
f_2, f_3$ be three real-valued functions over $\Z_q^*$. Then for
any $n\in\Z_q$, there exist $x, y, z\in\Z_q^*$ such that $n=x+y+z$
and
$$
f_1(x)+f_2(y)+f_3(z)\geqslant\frac{1}{\phi(q)}\sum_{a\in\Z_q^*}(f_1(a)+f_2(a)+f_3(a)),
$$
where $\phi$ is the Euler totient function.
\end{Thm}

The proof of Theorem \ref{t2} will be given in Section 2, and we
shall prove Theorem \ref{t1} in Section 3.

\section{Proof of Theorem \ref{t2}}
\setcounter{equation}{0} \setcounter{Thm}{0} \setcounter{Lem}{0}
\setcounter{Cor}{0}

Let
$$
K=\frac{1}{\phi(q)}\sum_{x\in\Z_q^*}(f_1(x)+f_2(x)+f_3(x)).
$$
We shall make an induction on the number of prime divisors of $q$.
First, assuming that Theorem \ref{t2} holds for two co-prime
integers $q_1$ and $q_2$, we claim that this theorem is also valid
for $q=q_1q_2$. Consider $\Z_q$ as $\Z_{q_1}\oplus\Z_{q_2}$ and
functions $g_1, g_2, g_3$ over $\Z_{q_1}^*$ by
$$
g_i(a)=\frac{1}{\phi(q_2)}\sum_{b\in\Z_{q_2}^*}f_i((a,b)).
$$
Thus for any $n=(n_1,n_2)\in\Z_{q_1}\oplus\Z_{q_2}$, by the
induction hypothesis, there exist $x_1, y_1, z_1\in\Z_{q_1}^*$
such that $n_1=x_1+y_1+z_1$ and
$$
g_1(x_1)+g_2(y_1)+g_3(z_1)\geqslant
\frac{1}{\phi(q_1)}\sum_{a\in\Z_{q_1}^*}(g_1(a)+g_2(a)+g_3(a))=K,
$$
i.e.,
$$
\sum_{b\in\Z_{q_2}^*}(f_1((x_1,b))+f_2((y_1,b))+f_3((z_1,b)))\geqslant
\phi(q_2)K.
$$
Define functions $h_1, h_2, h_3$ over $\Z_{q_2}^*$ by
$$
h_1(b)=f_1((x_1,b)),\ \ h_2(b)=f_2((y_1,b)),\ \text{and}\
h_3(b)=f_3((z_1,b)).
$$
Then applying the induction hypothesis again, there exist $x_2,
y_2, z_2\in\Z_{q_2}^*$ such that $n_2=x_2+y_2+z_2$ and
$$
h_1(x_2)+h_2(y_2)+h_3(z_2)\geqslant\frac{1}{\phi(q_2)}\sum_{b\in\Z_{q_2}^*}(h_1(b)+h_2(b)+h_3(b))\geqslant
K.
$$
This concludes the proof of our induction.

Thus we only need to prove Theorem \ref{t2} when $q$ is the power of
a prime. Assume that $q=p$ where $p\geqslant 5$ is a prime. Let $
S_i=\sum_{a\not=0}f_i(a)$ for $i=1,2,3$. Clearly
$S_1+S_2+S_3=(p-1)K$. Assume on the contrary that there exists some
$n\in\Z_p$ such that for any $x, y, z\in\Z_p^*$ with $x+y+z=n$,
\begin{equation}
\label{ie1} f_1(x)+f_2(y)+f_3(z)<K.
\end{equation}

We firstly consider the case $n=0$. Observe that
$$
\sum_{\substack{x,y,z\not=0\\
x+y+z=0}}f_1(x)=\sum_{x\not=0}f_1(x)\sum_{\substack{y,z\not=0\\
y+z=-x}}1=\sum_{x\not=0}f_1(x)\sum_{y\not=0,-x}1=(p-2)S_1.
$$
Similarly we have
$$
\sum_{\substack{x,y,z\not=0\\
x+y+z=0}}f_2(y)=(p-2)S_2\qquad\text{and}\qquad
\sum_{\substack{x,y,z\not=0\\
x+y+z=0}}f_3(z)=(p-2)S_3.
$$
Therefore
\begin{align*}
S_1+S_2+S_3=&\frac{1}{p-2}\sum_{\substack{x,y,z\not=0\\
x+y+z=0}}(f_1(x)+f_2(y)+f_3(z))\\
<&\frac{K}{p-2}\sum_{\substack{x,y,z\not=0\\
x+y+z=0}}1\\
=&(p-1)K,
\end{align*} which evidently leads to a
contradiction as desired.

Now suppose that $n\not=0$. Then for each $x\in\Z_p^*$
\begin{align*}
\sum_{\substack{y,z\not=0\\
y+z=n-x}}(f_2(y)+f_3(z))=&\sum_{y\not=0, n-x}(f_2(y)+f_3(n-x-y))\\
=&(S_2-f_2(n-x))+(S_3-f_3(n-x)),
\end{align*}
where we set $f_i(0)=0$. On the other hand, in view of
(\ref{ie1}),
$$
\sum_{\substack{y,z\not=0\\
y+z=n-x}}(f_2(y)+f_3(z))<(K-f_1(x))|\{(y,z):\, y,z\in\Z_p^*,
y+z=n-x\}|.
$$
Therefore
\begin{equation}
\label{ie2} S_2+S_3<(p-2)(K-f_1(x))+f_2(n-x)+f_3(n-x)
\end{equation}
for those $x\not=0, n$, and
\begin{equation}
\label{ie3} S_2+S_3<(p-1)(K-f_1(n)).
\end{equation}

Recalling that $S_1+S_2+S_3=(p-1)K$, we see that
$$
S_1>K+(p-2)f_1(x)-f_2(n-x)-f_3(n-x)
$$
provided that $x\not=0, n$. Summing the above inequality over all
$x\not=0, n$, we have
$$
(p-2)S_1>(p-2)K+(p-2)(S_1-f_1(n))-(S_2-f_2(n))-(S_3-f_3(n)),
$$
i.e.,
$$
S_2+S_3>(p-2)K-(p-2)f_1(n)+f_2(n)+f_3(n).
$$
Hence it follows from (\ref{ie2}) that
\begin{equation}
\label{ie4} (p-2)f_1(n)-f_2(n)-f_3(n)>(p-2)f_1(x)-f_2(n-x)-f_3(n-x),
\end{equation}
for any $x\not=0, n$. Symmetrically,
\begin{equation}
\label{ie5}
(p-2)f_2(n)-f_1(n)-f_3(n)>(p-2)f_2(x)-f_1(n-x)-f_3(n-x)
\end{equation}
and
\begin{equation}
\label{ie6}
(p-2)f_3(n)-f_1(n)-f_2(n)>(p-2)f_3(x)-f_1(n-x)-f_2(n-x).
\end{equation}
Computing $(p-3)\times (\ref{ie4})+(\ref{ie5})+(\ref{ie6})$, we
deduce that
\begin{align*}
(p-1)(p-4)f_1(n)>&(p-2)(p-3)f_1(x)-2f_1(n-x)+(p-2)f_2(x)\\
&-(p-2)f_2(n-x)+(p-2)f_3(x)-(p-2)f_3(n-x).
\end{align*}
Summing the above inequality over all $x\not=0,n$ again, then
$$
(p-2)(p-1)(p-4)f_1(n)>((p-2)(p-3)-2)(S_1-f_1(n)),
$$
i.e., $(p-1)f_1(n)>S_1$. Thus with the help of (\ref{ie3}), we
obtain a contradiction that
$$
S_1+S_2+S_3<S_1+(p-1)(K-f_1(n))<(p-1)K.
$$

Finally, suppose that $q=p^\alpha$ where $\alpha>1$. Define $g_1,
g_2, g_3$ over $\Z_p^*$ by $ g_i(x)=p^{1-\alpha}\sum_{a\equiv
x\pmod{p}}f_i(a)$. For any $n\in\Z_{p^\alpha}$, since Theorem
\ref{t2} holds for $p$, we know that there exist
$x_1,y_1,z_1\in\{1,2,\ldots,p-1\}$ such that $n\equiv
x_1+y_1+z_1\pmod{p}$ and
\begin{align*}
g_1(x_1)+g_2(y_1)+g_3(z_1)
\geqslant\frac{1}{p-1}\sum_{a\in\Z_{p}^*}(g_1(a)+g_2(a)+g_3(a))=K.
\end{align*}
Let $n'=(n-x_1-y_1-z_1)/p$, and define $h_1, h_2, h_3$ over
$\Z_{p^{\alpha-1}}$ by
$$
h_1(x)=f_1(x_1+xp),\ h_2(y)=f_2(y_1+yp)\text{ and
}h_3(z)=f_3(z_1+zp).
$$
It is easy to check that
$$
\sum_{\substack{x, y, z\in\Z_{p^{\alpha-1}}\\
x+y+z=n'}}(h_1(x)+h_2(y)+h_3(z))=p^{\alpha-1}\sum_{b\in\Z_{p^{\alpha-1}}}(h_1(b)+h_2(b)+h_3(b))
$$
and
$$
|\{(x,y,z):\, x, y, z\in\Z_{p^{\alpha-1}},
x+y+z=n'\}|=p^{2\alpha-2}.
$$
And we have
\begin{align*}
\sum_{b\in\Z_{p^{\alpha-1}}}(h_1(b)+h_2(b)+h_3(b))=&
\sum_{b\in\Z_{p^{\alpha-1}}}(f_1(x_1+bp)+f_2(y_1+bp)+f_3(z_1+bp))\\
=&p^{\alpha-1}(g_1(x_1)+g_2(y_1)+g_3(z_1))\\
\geqslant& p^{\alpha-1}K.
\end{align*}
Therefore there must exist $x_2, y_2,
z_2\in\Z_{p^{\alpha-1}}$ such that $n'=x_2+y_2+z_2$ and
$$
f_1(x_1+x_2p)+f_2(y_1+y_2p)+f_3(z_1+z_2p)=h_1(x_2)+h_2(y_2)+h_3(z_2)\geqslant
K.
$$
The proof is complete.\qed

\begin{Cor}
\label{c1} Let $q$ be a positive square-free odd integer. Suppose
that
$$
f_1,f_2,f_3:\ \Z_q^*\ \longrightarrow\ [0,1]
$$
satisfy that
$$
\sum_{a\in\Z_q^*}(f_1(a)+f_2(a)+f_3(a))>2\phi(q).
$$
Then for any $n\in\Z_q$, there exist $x,y,z\in\Z_q^*$ such that
$n=x+y+z$, $f_1(x)+f_2(y)+f_3(z)>5/3$ and
$f_1(x)f_2(y)f_3(z)\not=0$.
\end{Cor}
\begin{proof} In view of Theorem \ref{t2}, there is nothing to do
if $3\nmid q$. The case $q=3$ can be verified directly. For
example, supposing that $n=1$, we have
\begin{align*}
&\max\{f_1(1)+f_2(1)+f_3(2),f_1(1)+f_2(2)+f_3(1),f_1(2)+f_2(1)+f_3(1)\}\\
\geqslant&\frac{1}{3}(2f_1(1)+2f_2(1)+2f_3(1)+f_1(2)+f_2(2)+f_3(2))\\
>&\frac{8}{3}-\frac{1}{3}(f_1(2)+f_2(2)+f_3(2))\\
\geqslant&\frac{5}{3}.
\end{align*}
And if $f_1(1)=0$ (resp. $f_1(2)=0$), then $f_1(2)+f_2(1)+f_3(1)$
(resp. $f_1(1)+f_2(2)+f_3(1)$) is greater than
$4-f_2(2)-f_3(2)\geqslant 2$ (resp. $4-f_2(1)-f_3(2)\geqslant 2$).

Finally, assume that $q=3q'$ where $3\nmid q'$. By Theorem
\ref{t2}, for any $n=(n_1,n_2)\in\Z_{q'}\oplus\Z_3$ there exist
$x_1, y_1, z_1\in\Z_{q'}^*$ such that $n_1=x_1+y_1+z_1$ and
$$
\sum_{b\in\Z_3^*}(f_1((x_1,b))+f_2((y_1,b))+f_3((z_1,b)))>2\phi(3).
$$
It follows that there exist $x_2, y_2, z_2\in\Z_{3}^*$ such that
$n_2=x_2+y_2+z_2$,
$$
f_1((x_1,x_2))+f_2((y_1,y_2))+f_3((z_1,z_2))>\frac{5}{3}
$$
and
$$
f_1((x_1,x_2))f_2((y_1,y_2))f_3((z_1,z_2))\not=0.
$$
\end{proof}

\section{Proof of Theorem \ref{t1}}
\setcounter{equation}{0} \setcounter{Thm}{0} \setcounter{Lem}{0}
\setcounter{Cor}{0}

Our proof of Theorem \ref{t1} will follow that of Green in
\cite{Green05}, only with some slight modifications. Let
$$
\kappa=10^{-4}(\underline{d}_{\P}(P_1)+\underline{d}_{\P}(P_2)+\underline{d}_{\P}(P_3)-2),
$$
and let $\alpha_i=\underline{d}_{\P}(P_i)/(1+2\kappa)$. We may
assume that $n$ is sufficiently large so that
$$
|P_i\cap[1,2n/3]|\geqslant (1+\kappa)\alpha_i\frac{2n/3}{\log n}.
$$
Let $w=w(n)$ be a function tending sufficiently slowly to infinity
with $n$ (e.g., we may choose $w(n)=\floor{\frac{1}{4}\log\log
n}$), and let
$$
W=\prod_{\substack{p\in\P\\ p\leqslant w(n)}}p.
$$
Clearly $W\leqslant\log n$ and
\begin{align*}
\sum_{\substack{x\leqslant 2n/3\\ (x,W)=1}}\1_{P_i}(x)\log
x\geqslant& \sum_{n^{\frac{1}{1+\kappa/2}}\leqslant x\leqslant
2n/3}\1_{P_i}(x)\log
x\\
\geqslant&\frac{\log
n}{1+\kappa/2}\bigg(\frac{(1+\kappa)\alpha_i(2n/3)}{\log
n}-n^{\frac{1}{1+\kappa/2}}\bigg)\\
\geqslant&\frac{2}{3}\alpha_in,
\end{align*}
whenever $n$ is sufficiently large, where we set $\1_A(x)=1$ if
$x\in A$ and $0$ otherwise. Define
$$
f_i(b)=\max\bigg\{0,\frac{3\phi(W)}{2n}\sum_{\substack{x\leqslant 2n/3\\
x\equiv b\pmod{W}}}\1_{P_i}(x)\log x-3\kappa\bigg\}
$$
for $b\in\Z_W^*$. By the well-known Siegel-Walfisz theorem (cf.
\cite{Davenport00}), we know that $f_i(b)\in[0,1]$ if $n$ is
sufficiently large. Note that
\begin{align*}
&\sum_{b\in\Z_W^*}(f_1(b)+f_2(b)+f_3(b))\\
\geqslant&\frac{3\phi(W)}{2n}\sum_{b\in\Z_W^*}\sum_{\substack{x\leqslant 2n/3\\
x\equiv b\pmod{W}}}(\1_{P_1}(x)+\1_{P_2}(x)+\1_{P_3}(x))\log
x-9\kappa\phi(W)\\
\geqslant&(\alpha_1+\alpha_2+\alpha_3-9\kappa)\phi(W)\\
>&2\phi(W).
\end{align*}
In view of Corollary \ref{c1}, there exist $b_1, b_2,
b_3\in\Z_W^*$ such that $n\equiv b_1+b_2+b_3\pmod{W}$,
$f_1(b_1)+f_2(b_2)+f_3(b_3)>5/3$ and $f_i(b_i)>0$. And without
loss of generality, we may assume that $1\leqslant b_1, b_2,
b_3<W$.

Let $N$ be a prime in the interval
$[(1+\kappa)n/W,(1+2\kappa)n/W]$. Thanks to the prime number
theorem, such $N$ always exists for sufficiently large $n$.
Following our discussions above, let $n'=(n-b_1-b_2-b_3)/W$ and
let
$$
A_i=\{x:\, Wx+b_i\in P_i\cap[1,2n/3]\}.
$$
It suffices to show that $n'\in A_1+A_2+A_3$. Let
$$
\alpha_i'=\sum_{x}\1_{A_i}(x)\lambda_{b_i,W,N}(x),
$$
where
$$
\lambda_{b,W,N}(x)=\begin{cases} \phi(W)\log(Wx+b)/WN&\text{ if
}x\leqslant N\text{ and }Wx+b\text{ is prime},\\
0&\text{otherwise}.
\end{cases}
$$
Note that
$$
\alpha_i'=\frac{\phi(W)}{WN}\sum_{x}\1_{A_i}(x)\log(Wx+b_i)\geqslant
\frac{2(f_i(b_i)+3\kappa)}{3(1+2\kappa)},
$$
since $A_i\subseteq[0,N]$ and $f_i(b_i)>0$. Then we have
$$
\alpha_1', \alpha_2', \alpha_3'\geqslant
\frac{2\kappa}{1+2\kappa}\geqslant\kappa
$$
and
\begin{align*}
\alpha_1'+\alpha_2'+\alpha_3'\geqslant&\frac{2}{3(1+2\kappa)}(f_1(b_1)+f_2(b_2)+f_3(b_3)+9\kappa)\\
\geqslant&\frac{10}{9}+3\kappa.
\end{align*}

Below we consider $A_1, A_2, A_3$ as the subsets of $\Z_N$. Since
$A_1, A_2, A_3\subseteq[0,2n/3W]$ and $N\geqslant n/W+3$, there
exist no $x_i\in A_i$ such that $x_1+x_2+x_3=n'+N$ in $\Z$.
Therefore $n'\in A_1+A_2+A_3$ in $\Z_N$ implies that $n'\in
A_1+A_2+A_3$ in $\Z$. Let $\mu_i(x)=\lambda_{b_i,W,N}(x)$ and
$a_i(x)=\1_{A_i}(x)\mu_i(x)$. For an arbitrary complex-valued
function $f$ over $\Z_N$, define $\tilde{f}$ over $\Z_N$ by
$$
\tilde{f}(r)=\sum_{x\in\Z_N}f(x)e(-xr/N),
$$
where $e(x)=e^{2\pi\sqrt{-1}x}$. Also, for functions $f, g$ over
$\Z_N$, define
$$
(f*g)(x)=\sum_{y\in\Z_N}f(y)g(x-y).
$$
It is easy to check that $(f*g)\,\tilde{}=\tilde{f}\tilde{g}$.
Suppose that $\delta, \epsilon>0$ are two real numbers which will
be chosen later. Let
$$
R_i=\{r\in\Z_N:\, |\tilde{a_i}(r)|\geqslant\delta\}
$$
and
$$
B_i=\{x\in\Z_N:\, \|xr/N\|\leqslant\epsilon\text{ for all }r\in
R_i\},
$$
where $ \|x\|=\min_{z\in\Z}|x-z|$. Also let $\beta_i=\1_{B_i}/|B_i|$
and $a_i'=a_i*\beta_i*\beta_i$.
\begin{Lem}
\label{l1}
$$
\bigg|\sum_{\substack{x,y,z\in\Z_N\\
x+y+z=n'}}a_1'(x)a_2'(y)a_3'(z)-\sum_{\substack{x,y,z\in\Z_N\\
x+y+z=n'}}a_1(x)a_2(y)a_3(z)\bigg|\leqslant\frac{C_1}{N}(\epsilon^2\delta^{-5/2}+\delta^{1/4}).
$$
\end{Lem}
\begin{proof}
It is not difficult to see that
$$
\sum_{\substack{x,y,z\in\Z_N\\
x+y+z=n'}}a_1(x)a_2(y)a_3(z)=N^{-1}\sum_{r\in\Z_N}e(n'r/N)\tilde{a_1}(r)\tilde{a_2}(r)\tilde{a_3}(r).
$$
Thus
\begin{align*}
&\bigg|\sum_{\substack{x,y,z\in\Z_N\\
x+y+z=n'}}a_1'(x)a_2'(y)a_3'(z)-\sum_{\substack{x,y,z\in\Z_N\\
x+y+z=n'}}a_1(x)a_2(y)a_3(z)\bigg|\\
=&N^{-1}\bigg|\sum_{r\in\Z_N}e(n'r/N)\tilde{a_1}(r)\tilde{a_2}(r)\tilde{a_3}(r)(1-\tilde{\beta_1}^2(r)\tilde{\beta_2}^2(r)\tilde{\beta_3}^2(r))\bigg|.
\end{align*}
From the proofs of Lemma 6.7 and Proposition 6.4 in
\cite{Green05}, we know that $|R_i|\leqslant C_1'\delta^{-5/2}$
for some absolute constant $C_1'$, and if $r\in R=R_1\cap R_2\cap
R_3$, then
$$
|1-\tilde{\beta_1}^2(r)\tilde{\beta_2}^2(r)\tilde{\beta_3}^2(r)|\leqslant
2^{12}\epsilon^2.
$$
Therefore
\begin{align*}
&\bigg|\sum_{r\in
R}e(n'r/N)\tilde{a_1}(r)\tilde{a_2}(r)\tilde{a_3}(r)(1-\tilde{\beta_1}^2(r)\tilde{\beta_2}^2(r)\tilde{\beta_3}^2(r))\bigg|\\
\leqslant&2^{12}\epsilon^2\sum_{r\in
R}\bigg|e(n'r/N)\tilde{a_1}(r)\tilde{a_2}(r)\tilde{a_3}(r)\bigg|\\\leqslant&
2^{12}\epsilon^2|R|\\\leqslant& 2^{12}C_1'\epsilon^2\delta^{-5/2}
\end{align*}
by noting that
$|\tilde{a}_i(r)|\leqslant\sum_{x\in\Z_N}a_i(x)\leqslant 1$. And
since $|\tilde{\beta}_i(r)|\leqslant\sum_{x\in\Z_N}\beta_i(x)=1$,
with the help of the H\"older inequality,
\begin{align*}
&\bigg|\sum_{r\not\in
R}e(n'r/N)\tilde{a_1}(r)\tilde{a_2}(r)\tilde{a_3}(r)(1-\tilde{\beta_1}^2(r)\tilde{\beta_2}^2(r)\tilde{\beta_3}^2(r))\bigg|\\
\leqslant&2\sup_{r\not\in
R}|\tilde{a_1}(r)\tilde{a_2}(r)\tilde{a_3}(r)|^{1/4}\sum_{r\not\in
R}|\tilde{a_1}(r)|^{3/4}|\tilde{a_2}(r)|^{3/4}|\tilde{a_3}(r)|^{3/4}\\
\leqslant&2\delta^{1/4}\bigg(\sum_{r\not\in
R}|\tilde{a_1}(r)|^{9/4}\bigg)^{1/3}\bigg(\sum_{r\not\in
R}|\tilde{a_2}(r)|^{9/4}\bigg)^{1/3}\bigg(\sum_{r\not\in
R}|\tilde{a_3}(r)|^{9/4}\bigg)^{1/3}\\
\leqslant&C_1''\delta^{1/4},
\end{align*}
where we apply Lemma 6.6 in \cite{Green05} with $p=9/4$. This
concludes our proof.
\end{proof}
Now we shall give a lower bound only depending on $\kappa$ for
$\sum_{x+y+z=n'}a_1'(x)a_2'(y)a_3'(z)$.
\begin{Lem}
\label{l2} Suppose that $\epsilon^{|R_i|}\geqslant C\log\log w/w$.
Then for each $x\in\Z_N$
$$
|a_i'(x)|\leqslant (1+2C^{-1})/N.
$$
\end{Lem}
\begin{proof} The proof is same as Lemma 6.3 in \cite{Green05}, so
we omit the details here.
\end{proof}

In \cite{Varnavides59}, Varnavides showed that if $A$ is a subset
of $\Z_N$ with $|A|\geqslant \theta N$, then $A$ contains at least
$c(\theta)N^2$ non-trivial 3APs whenever $N$ is sufficiently
large, where $c(\theta)$ is a constant only depending on $\theta$.
Varnavides' argument was used by Green in the proof of his Lemma
6.8 \cite{Green05}. Here we also need an analogue of Varnavides'
result for sumsets. For non-empty subsets $X_1, X_2, \ldots, X_k$
of $\Z_N$, define
$$
\nu_{X_1,X_2,\ldots,X_k}(n)=|\{(x_1,x_2,\ldots,x_k):\, x_i\in X_i,
n=x_1+x_2+\cdots+x_k\}|.
$$
In particular, we set $\nu_{X_1}(n)=\1_{X_1}(n)$.
\begin{Lem}
\label{l3} Suppose that $k\geqslant 2$ and
$0<\theta_1,\ldots,\theta_k\leqslant 1$ with
$\theta_1+\cdots+\theta_k>1$. Let
$$
\theta=\min\{\theta_1,\ldots,\theta_k,(\theta_1+\cdots+\theta_k-1)/(3k-5)\}.
$$
Suppose that $N$ is a prime greater than $2\theta^{-2}$, and
$X_1,\ldots,X_k$ are subsets of $\Z_N$ with
$|X_i|\geqslant\theta_i N$. Then for any $n\in\Z_N$, we have
$\nu_{X_1,X_2,\ldots,X_k}(n)\geqslant\theta^{2k-3}N^{k-1}$.
\end{Lem}
\begin{proof}
When $k=2$, we have
$$
\nu_{X_1,X_2}(n)=|X_1\cap(n-X_2)|\geq
|X_1|+|X_2|-|X_1\cup(n-X_2)|\geq|X_1|+|X_2|-N.
$$
Below we assume that $k\geqslant 3$ and the assertion holds for
the smaller values of $k$.

Suppose that $A, B$ are two non-empty subsets of $\Z_N$. Let
$$
S_t(A,B)=\{x\in\Z_N:\, \nu_{A,B}(x)\geqslant t\}.
$$
A result of Pollard \cite{Pollard74, Nathanson96} asserts that for
any $1\leqslant t\leqslant\min\{|A|,|B|\}$
$$
\sum_{i=1}^t|S_i(A,B)|\geqslant\min\{tN,t(|A|+|B|-t)\}.
$$
(The case $t=1$ is the well-known Cauchy-Davenport theorem.)

Without loss of generality, we suppose that
$\theta_1\geqslant\theta_2\geqslant\cdots\geqslant\theta_k$. If
$\theta_1+\theta_2>1+\theta$, then
\begin{align*}
\nu_{X_1,X_2,\ldots,X_k}(n)=&\sum_{x\in
X_3+\cdots+X_k}\nu_{X_1,X_2}(n-x)
\nu_{X_3,\cdots,X_k}(x)\\\geqslant&(|X_1|+|X_2|-N)\sum_{x\in
X_3+\cdots+X_k}\nu_{X_3,\cdots,X_k}(x)\\
=&(|X_1|+|X_2|-N)|X_3|\cdots|X_k|\\
\\\geqslant&(\theta_1+\theta_2-1)\theta_3\cdots\theta_kN^{k-1}.
\end{align*}

Now we may assume that $\theta_1+\theta_2\leqslant 1+\theta$. Let
$\ceil{t}=\min\{n:n\in \mathbb{Z},n\geq t\}$, by Pollard's
theorem, we have
$$
\sum_{i=1}^{\ceil{\theta N}}|S_i(X_1,X_2)|\geqslant\ceil{\theta
N}(\theta_1N+\theta_2N-\ceil{\theta N}).
$$
It follows that
$$
\sum_{i=\ceil{\theta^2 N}}^{\ceil{\theta
N}}|S_i(X_1,X_2)|\geqslant\ceil{\theta
N}(\theta_1N+\theta_2N-\ceil{\theta N})-\ceil{\theta^2 N}N.
$$
Hence by noting that
$$
\frac{\ceil{\theta N}}{N}+ \frac{\ceil{\theta^2 N}}{\ceil{\theta N}}
\leqslant 2\theta+\frac{1}{N}+\frac{1}{\theta N} \leqslant 3\theta,
$$
we have
\begin{align*}
|S_{\ceil{\theta^2 N}}(X_1,X_2)|\geqslant&\frac{\ceil{\theta
N}(\theta_1N+\theta_2N-\ceil{\theta N})-\ceil{\theta^2
N}N}{\ceil{\theta N}-\ceil{\theta^2 N}+1}\\
\geqslant& \theta_1N+\theta_2N-3\theta N.
\end{align*}
Let $Y=S_{\ceil{\theta^2 N}}(X_1,X_2)$. Clearly $|Y|\geqslant\theta
N$ since
$$
\theta\leqslant
\frac{\theta_1+\theta_2+\cdots+\theta_k-1}{3k-5}\leqslant\frac{k-1}{2(3k-5)}(\theta_1+\theta_2).
$$
Then by the induction hypothesis on $k$,
\begin{align*}
\nu_{X_1,X_2,\ldots,X_k}(n)\geqslant&
\sum_{x\in Y}\nu_{X_1,X_2}(x)\nu_{X_3,\ldots,X_k}(n-x)\\
\geqslant&\inf_{x\in Y}\nu_{X_1,X_2}(x)\sum_{x\in
Y}\nu_{X_3,\ldots,X_k}(n-x)\\
=&\nu_{Y,X_3,\ldots,X_k}(n)\inf_{x\in Y}\nu_{X_1,X_2}(x)\\
\geqslant&\theta_*^{2k-5}N^{k-2}\theta^2N,
\end{align*}
where
$$
\theta_*=\min\{\theta_1+\theta_2-3\theta,\theta_3,\ldots,\theta_k,(\theta_1+\cdots+\theta_k-3\theta-1)/(3k-8)\}\geqslant\theta.
$$
\end{proof}
\begin{Lem}
\label{l4}
$$
\sum_{\substack{x,y,z\in\Z_N\\
x+y+z=n'}}a_1'(x)a_2'(y)a_3'(z)\geqslant \frac{\kappa^9}{8N}.
$$
\end{Lem}
\begin{proof}
Let $A_i'=\{x\in\Z_N:\, a_i'(x)\geqslant\alpha_i'\kappa/N\}$.
Applying Lemma \ref{l2} with $C=2/\kappa$, we have
$$
\alpha_i'=\sum_{x\in\Z_N}a_i(x)=\sum_{x\in\Z_N}a_i'(x)\leqslant
\frac{1+\kappa}{N}|A_i'|+\frac{\alpha_i'\kappa}{N}(N-|A_i'|),
$$
whence
$$
|A_i'|\geqslant\frac{\alpha_i'(1-\kappa)}{1+\kappa}N.
$$
Observe that $\alpha_i'(1-\kappa)/(1+\kappa)\geqslant \kappa/2$
and
$$
\sum_{i=1}^3\frac{\alpha_i'(1-\kappa)}{1+\kappa}=\frac{1-\kappa}{1+\kappa}(\alpha_1'+\alpha_2'+\alpha_3')\geqslant
\frac{10}{9}+\frac{\kappa}{2}.
$$
Then with the help of Lemma \ref{l3},
$$
\nu_{A_1',A_2',A_3'}(n')\geqslant \frac{\kappa^3}{8}N^2.
$$
It follows that
$$
\sum_{\substack{x,y,z\in\Z_N\\
x+y+z=n'}}a_1'(x)a_2'(y)a_3'(z)\geqslant\sum_{\substack{x\in A_1',y\in A_2',z\in A_3'\\
x+y+z=n'}}a_1'(x)a_2'(y)a_3'(z)\geqslant
\frac{\kappa^6}{8}\alpha_1'\alpha_2'\alpha_3'N^{-1}.
$$
\end{proof}

Now combining Lemmas \ref{l1} and \ref{l4}, we obtain that
$$
N\sum_{\substack{x,y,z\in\Z_N\\
x+y+z=n'}}a_1(x)a_2(y)a_3(z)+C_1(\epsilon^2\delta^{-5/2}+\delta^{1/4})\geqslant\frac{\kappa^9}{8}.
$$
By the final arguments in \cite{Green05}, we know that under the
condition in Lemma \ref{l1}, we may choose $\delta$ and $\epsilon$
such that both $\epsilon^2\delta^{-5/2}$ and $\delta^{1/4}$ tend
to $0$, whenever $N$ is sufficiently large. Thus for sufficiently
large $n$,
$$
N\sum_{x+y+z=n'}a_1(x)a_2(y)a_3(z)\geqslant\frac{\kappa^9}{9}>0.
$$
\qed

\section{Further Remarks}
\setcounter{equation}{0} \setcounter{Thm}{0} \setcounter{Lem}{0}
\setcounter{Cor}{0}

Maybe the most famous unsolved conjecture in number theory is the
binary Goldbach problem, which says that every even integer greater
than 2 is the sum of two primes. The well-known result of Chen
\cite{Chen73} asserts that every sufficiently large even integer can
be represented as the sum of a prime and an integer which is a prime
or the product of two primes. However, it seems that a similar
extension of above result will fail for the binary Goldbach
conjecture. For any $\epsilon>0$, there exists a sufficiently large
$w$ such that
$$
\prod_{\substack{p\text{ prime}\\ 3\leq p\leq
w}}\bigg(\frac{p-2}{p-1}\bigg)<\epsilon.
$$
Let $P_1=\P\cap(w,\infty)$ and
$$
P_2=\{x\in\P:\,x\equiv 1\pmod{p}\text{ for an odd prime }p\leq
w\}.
$$
Clearly
$$
\underline{d}_{\P}(P_1)+\underline{d}_{\P}(P_2)=1+1-\prod_{\substack{p\text{ prime}\\
3\leq p\leq w}}\bigg(\frac{p-2}{p-1}\bigg)>2-\epsilon.
$$
But $Wk+1\not\in P_1+P_2$ for each odd integers $k$, where
$$
W=\prod_{\substack{p\text{ prime}\\ 3\leq p\leq w}}p.
$$

In fact, we can construct two sets $P_1, P_2$ of primes with
$\underline{d}_{\P}(P_1)=\underline{d}_{\P}(P_2)=1$ such that
there exist infinitely many positive even integers not contained
in $P_1+P_2$. Let $N_k=2\lfloor e^{k\sqrt{\log k}}\rfloor$ and
$n_k=N_{k+1}+N_k+2$. Let ${\mathcal A}_k=\{p\in\P:\, n_k-p\in\P\}
$ and
$$
{\mathcal B}_k=\{p\in\P:\, N_k+2N_{k-1}\leqslant p\leqslant
N_{k+1}, n_k-p\not\in\P\}.
$$
Set $P_1=P_2=\bigcup_{k=1}^{\infty}{\mathcal B}_k$. With the help
of Selberg's sieve method, we know that
$$
|{\mathcal A}_k|\ll \frac{n_k}{(\log n_k)^2}\prod_{\substack{p\text{
prime}\\ p\mid n_{k}}}\bigg(1+\frac{1}{p}\bigg).
$$
Define $z(n)=\max\{z\in\N:\, \prod_{\substack{p\text{ prime}\\
p\leqslant z}}p\leqslant n\}$. By the prime number theorem, $z(n)\ll
\log n$. Hence by the Mertens theorem,
$$
\prod_{\substack{p\text{ prime}\\ p\mid
n}}\bigg(1+\frac{1}{p}\bigg)\leqslant\prod_{\substack{p\text{ prime}\\
p\leqslant z(n)}}\bigg(1+\frac{1}{p}\bigg)\ll\log\log n.
$$
So $|{\mathcal A}_k|\ll N_{k+1}\log\log N_{k+1}/(\log N_{k+1})^2$.
It is not difficult to verify that
$$
\frac{N_{k}}{\log N_{k}}=o\bigg(\frac{N_{k+1}}{\log
N_{k+1}}\bigg)\text{\quad and\quad} \frac{N_{k+1}(\log\log
N_{k+1})^2}{(\log N_{k+1})^2}=o\bigg(\frac{N_{k}}{\log N_k}\bigg).
$$
So by the prime number theorem, for $x\in(N_k,N_k+2N_{k-1})$ we
have
\begin{align*}
|P_1\cap[1,x]|\geqslant&|{\mathcal B}_{k-1}|\geqslant|\P\cap[N_{k-1}+2N_{k-2},N_k]|-|{\mathcal A}_{k-1}|\\
\geqslant&(1+o(1))\bigg(\frac{N_{k}}{\log
N_{k}}-\frac{N_{k-1}+2N_{k-2}}{\log(N_{k-1}+2N_{k-2})}\bigg)-o\bigg(\frac{N_{k}}{\log
N_{k}}\bigg)\\
=&(1+o(1))\frac{x}{\log x}.
\end{align*}
And for $x\in[N_k+2N_{k-1},N_{k+1}]$,
\begin{align*}
|P_1\cap[1,x]|\geqslant& |{\mathcal B}_{k}\cap[1,x]|+|{\mathcal B}_{k-1}|\\
\geqslant&|\P\cap[N_{k}+2N_{k-1},x]|-|{\mathcal A}_k|+|{\mathcal B}_{k-1}|\\
\geqslant&(1+o(1))\bigg(\frac{x}{\log
x}-\frac{N_k+2N_{k-1}}{\log(N_k+2N_{k-1})}+\frac{N_{k}}{\log
N_k}\bigg)\\
=&(1+o(1))\frac{x}{\log x}.
\end{align*}
It follows that
$\underline{d}_{\P}(P_1)=\underline{d}_{\P}(P_2)=1$. But now
$n_k\not\in P_1+P_2$, since
$$
n_k\not\in\bigg(\bigcup_{j\not=k}[N_{j}+2N_{j-1},N_{j+1}]\bigg)+\bigg(\bigcup_{j}[N_{j}+2N_{j-1},N_{j+1}]\bigg)
$$
and $n_k\not\in {\mathcal B}_k+{\mathcal B}_k$.

Moreover, we mention that $\underline{d}$ can't be replaced by
$\overline{d}$ in Theorem \ref{t1}. Let $N_k=2\lfloor
e^{k\sqrt{\log k}}\rfloor$ and
$$
{\mathcal A}_k=\{n:\,2\mid n, N_{k+1}+N_{k}+2\leqslant n\leqslant
N_{k+1}+N_{k}+2\lfloor\log\log N_{k+1}\rfloor\}.
$$
Let
$$
{\mathcal B}_k=\{p\in\P:\, n-p\in\P\text{ for some }n\in {\mathcal
A}_k\}
$$
and
$$
{\mathcal C}_k=\{p\in\P:\, N_k+2N_{k-1}\leqslant p\leqslant
N_{k+1}, n-p\not\in\P\text{ for every }n\in {\mathcal A}_k\}.
$$
Then
$$
|{\mathcal B}_k|=O\bigg(\frac{N_{k+1}(\log\log N_{k+1})^2}{(\log
N_{k+1})^2}\bigg)=o\bigg(\frac{N_k}{\log N_k}\bigg).
$$
Let $P_1=P_2=\bigcup_{k=1}^\infty {\mathcal C}_k$. Similarly as
above, we also have
$\underline{d}_{\P}(P_1)=\underline{d}_{\P}(P_2)=1$ and $n\not\in
P_1+P_2$ for any $n\in {\mathcal A}_k$. Let $M_1=2$ and
$M_{l+1}=e^{e^{M_{l}}}$. Let
$$
P_3=\P\cap\bigg(\bigcup_{l=1}^{\infty}[M_{3l},M_{3l+1}]\bigg).
$$
Evidently $\overline{d}_{\P}(P_3)=1$. And for sufficiently large
$l$, there always exists $k$ such that $M_{3l+2}<
N_{k+1}<M_{3l+3}/2$. Let $n_k=N_{k+1}+N_k+2\lfloor\log\log
N_{k+1}\rfloor-1$. Assume that $n_k=p_1+p_2+p_3$ where $p_i\in
P_i$. Then we must have $p_3\leqslant M_{3l+1}$ since
$n_k\leqslant 2N_{k+1}<M_{3l+3}$. Hence $n_k-p_3\in {\mathcal
A}_k$ by noting $M_{3l+1}\leqslant \log\log M_{3l+2}<\log\log
N_{k+1}$. This leads to a contradiction since ${\mathcal
A}_k\cap(P_1+P_2)=\emptyset$.

\end{document}